\documentclass{article}
\usepackage{amsthm,amsmath, amssymb, fullpage, hyperref, tikz, caption}

\newenvironment{clmproof}[1]{\par\noindent\underline{Proof.}\space#1}{\leavevmode\unskip\penalty9999\hbox{}\nobreak\hfill\quad\hbox{$\diamondsuit$}\smallskip}

\RequirePackage{marginnote,hyperref}
\newcommand{\aside}[1]{\marginnote{\scriptsize{#1}}[0cm]}

\newcommand\Emph[1]{\emph{#1}\aside{#1}}

\renewcommand\ge\geqslant
\renewcommand\le\leqslant

\newcommand{\udot}{\mathbin{\mathaccent\cdot\cup}}

\theoremstyle{plain}
\newtheorem{thm}{Theorem}[section]

\newtheorem{lem}[thm]{Lemma}

\newtheorem{cor}[thm]{Corollary}
\newtheorem{clm}{Claim}
\newtheorem{rem}{Remark}

\newtheorem*{lemA}{Lemma A}
\newtheorem*{clique-corollary}{Corollary~\ref{clique-cor}}

\theoremstyle{definition}
\newtheorem{defn}{Definition}
\newtheorem*{clique-lem}{Big Clique Lemma}
\newcommand\C{\mathcal{C}}
\newcommand\dist{\mbox{dist}}

\title{Bounding Clique Size in Squares of Planar Graphs}
\author{Daniel W. Cranston\thanks{Virginia Commonwealth University, Department
of Computer Science, Richmond, VA. \texttt{dcranston@vcu.edu}}}

\begin{document}
\maketitle
\abstract{Wegner conjectured that if $G$ is a planar graph with maximum degree
$\Delta\ge 8$, then $\chi(G^2)\le \left\lfloor \frac32\Delta\right\rfloor +1$.
This problem has received much attention, but remains open for all $\Delta\ge 8$.  
Here we prove an analogous bound on $\omega(G^2)$: If $G$ is a plane graph with
$\Delta(G)\ge 36$, then $\omega(G^2)\le \lfloor\frac32\Delta(G)\rfloor+1$.
In fact, this is a corollary of the following lemma, which is our main result.
If $G$ is a plane graph with $\Delta(G)\ge 19$ and $S$ is a maximal clique in
$G^2$ with $|S|\ge \Delta(G)+20$, then there exist $x,y,z\in V(G)$ such that
$S=\{w:|N[w]\cap\{x,y,z\}|\ge 2\}$.}

\section{Introduction}
The square $G^2$ of a graph $G$ is formed from $G$ by adding an edge $vw$ for
each pair $v,w$ such that $\dist_G(v,w)=2$.  Easy counting arguments give
$\Delta(G)^2+1 \ge \chi(G^2)\ge \omega(G^2)\ge \Delta(G)+1$.  The first
inequality can hold with equality.  And in fact there exists an infinite
sequence of graphs $G_i$, with $\Delta(G_i)\to \infty$, such that
$\lim_{i\to\infty}\chi(G_i^2)/(\Delta(G_i)^2)=1$.  However, for planar graphs we can
do much better.

Wegner conjectured that if $G$ is planar and $\Delta(G)\ge 8$, then
$\chi(G^2)\le \lfloor\frac32\Delta(G)\rfloor+1$.  He also constructed
examples showing this bound would be best possible.  Form $G_{2s}$ as follows.
Beginning with vertices $x,y,z$, add $s$ common neighbors of each pair in
$\{x,y,z\}$, and then contract an edge incident to one common neighbor 
of $x$ and $y$ and to one common neighbor of $x$ and $z$.  The resulting graph
$G_{2s}$ has maximum degree $2s$, has order $3s+1$, and has diameter 2.
Thus, $G_{2s}^2=K_{3s+1}$. (The case when $\Delta$ is odd is similar.)

A graph $G$ is \emph{$k$-degenerate} if every subgraph $H$ of $G$ has a vertex
$v$ such that $d_H(v)\le k$.  If $G$ is $k$-degenerate, then repeatedly deleting
a vertex of degree at most $k$ gives a \emph{$k$-degeneracy order}.
Coloring greedily in the reverse of this order shows that
$\chi(G)\le k+1$.  If $G$ is $k$-degenerate, then this is witnessed by some
$k$-degeneracy order $\sigma$.  And the same order $\sigma$ witnesses that $G^2$
has degeneracy at most $\Delta(G)(2k-1)$.  Since every planar graph $G$ is
5-degenerate, we also get that $G^2$ has degeneracy at most $9\Delta(G)$, so
$\chi(G^2)\le 9\Delta(G)+1$.
Numerous papers have improved this bound~\cite{ AH, BBGvdH2, jonas, MS, vdHM03, Wong}.  
The current strongest result in this direction is $\chi(G^2)=\frac32(\Delta(G)+o(1))$.
 This was proved be Havet,
van den Heuvel, McDiarmid, and Reed~\cite{HvdHMR}.  Later, the proof was
simplified by Amini, Esperet, and van den Heuvel~\cite{AEvdH}.  For a recent
survey on coloring squares of graphs, particularly planar graphs, see~\cite{squares-survey}.

Despite this significant progress, Wegner's conjecture remains open for all
$\Delta\ge 8$.  This motivates interest in bounding $\omega(G^2)$ for all planar
graphs $G$.  Amini et al.~\cite{AEvdH} showed that for every surface $S$ there
exists a constant $C_S$ such that $\omega(G^2)\le\frac32\Delta(G)+C_S$.
When $S$ is the plane, they proved if $\Delta(G)\ge 11616$, then $\omega(G^2)\le
\frac32\Delta(G)+76$.  We prove a sharp upper bound on $\omega(G^2)$ when
$G$ is planar.  More
generally, we characterize all big cliques in squares of planar
graphs.  The following is our main result.  

\begin{clique-lem}
Fix an integer $D\ge 19$ and a plane graph $G$ with $\Delta(G)\le D$. If
$S\subseteq V(G)$ such that $S$ is a maximal clique in $G^2$ and $|S|\ge
D+20$, then $G$ contains vertices $x,y,z$ such that $S=\{v\in V(G):
|N[v]\cap\{x,y,z\}|\ge 2\}$.
\end{clique-lem}

Analogous to the bound Wegner conjectured on $\chi(G^2)$,
the Big Clique Lemma implies the following bound on $\omega(G^2)$.%
\footnote{It was reported in~\cite{AEvdH}, that Corollary~\ref{clique-cor} 
was proved by Cohen and van den Heuvel.  However, Cohen and van den Heuvel 
confirmed that they have not written a proof and have no plans to do so.}
This result is sharp, as witnessed by $G_{2s}$, constructed above.

\begin{cor}
\label{clique-cor}
If $G$ is a plane graph with $\Delta(G)\ge 36$, then $\omega(G^2)\le
\left\lfloor\frac32\Delta(G)\right\rfloor+1$.
\end{cor}

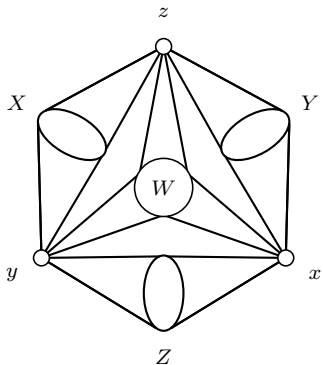
\begin{figure}[!h]
\centering
\begin{tikzpicture}[thick, scale=.75]
\tikzstyle{uStyle}=[shape = circle, minimum size = 6pt, inner sep = 0pt,
outer sep = 0pt, draw, fill=white, semithick]
\tikzstyle{sStyle}=[shape = rectangle, minimum size = 4.5pt, inner sep = 0pt,
outer sep = 0pt, draw, fill=white, semithick]
\tikzstyle{lStyle}=[shape = circle, minimum size = 4.5pt, inner sep = 0pt,
outer sep = 0pt, draw=none, fill=none]
\tikzset{every node/.style=uStyle}
\def\rad{2.5cm}
\def\off{.6cm}

\foreach \ang/\name in {90/z, 210/y, 330/x}
\draw (\ang:\rad) node[uStyle] (\name) {} ++ (\ang:\off) node[lStyle] {\footnotesize{$\name$}};

\foreach \ang/\name in {150/X, 270/Z, 30/Y}
{
\draw[rotate=\ang] (0:.75*\rad) ellipse (19pt and 10pt) {}; 
\draw (\ang:1.2*\rad) node[lStyle] {\footnotesize{$\name$}};
}

\foreach \name/\ang in {y/153, y/267, x/273, x/27, z/33, z/147}
{
\draw (\name) -- (\ang:\rad) (\name) -- (\ang:.49*\rad);
\draw (\name) -- (\ang:\rad) (\name) -- (\ang:.2*\rad);
}

\draw (0,0) node[circle, minimum size=22pt] {\footnotesize{$W$}};

\end{tikzpicture}
\caption{A partition of $S$ in the proof of Corollary~\ref{clique-cor}.\label{cor-fig}}
\end{figure}

\begin{proof}
Suppose the corollary is false.  Let $G$ be a counterexample, let
$D:=\Delta(G)$, and let $S$ be a maximum clique in $G^2$; note that $|S|\ge
\lfloor\frac32D\rfloor+2\ge D+20$, since $D\ge 36$.  By the Big Clique Lemma,
there exist vertices $x,y,z\in V(G)$ such that $S=\{v\in V(G):
|N[v]\cap\{x,y,z\}|\ge 2\}$.
Let $T:=S\cap\{x,y,z\}$, let $W:=N(x)\cap N(y)\cap N(z)$,
let $X:=N(y)\cap N(z)\setminus (T\cup W)$, let $Y:=N(x)\cap N(z)\setminus (T\cup
W)$, and let $Z:=N(x)\cap N(y)\setminus (T\cup W)$; see Figure~\ref{cor-fig}. 
So $|S|= |T|+|W|+|X|+|Y|+|Z|$.
Now $|Y|+|Z|+|W|+d_T(x)\le d(x)\le D$ and $|X|+|Z|+|W|+d_T(y)\le d(y)\le D$ and
$|X|+|Y|+|W|+d_T(z)\le d(z)\le D$.
Summing these 3 inequalities gives $2(|X|+|Y|+|Z|)+3|W|+\sum_{w\in T}d_T(w) \le 3D$.  
If $T\ne \emptyset$, then $|T|\in\{2,3\}$ and $G[T]$ is connected.  Thus, by
the Handshaking Lemma, the final sum on the left is at least $2(|T|-1)$.  So we
get $|S|=|X|+|Y|+|Z|+|T|+|W|\le \left\lfloor\frac32D\right\rfloor+1$.  
\end{proof}

Our definitions are mainly standard, but we highlight a few.
All of our graphs are finite and simple (with no loops).
For a graph $G$, we write $\Delta(G)$ and $\omega(G)$ to denote its maximum
degree and clique number.  When $G$ is clear from context, we may simply
write $\Delta$ or $\omega$.  A $k$-vertex (resp.~$k^+$-vertex or $k^-$-vertex)
is a vertex of degree exactly $k$ (resp.~at least $k$ or at most $k$).
We write $\dist_G(v,w)$ to denote the distance in $G$ between vertices $v$ and
$w$.  A \Emph{plane graph} is a fixed planar embedding of a planar graph.
We will also need the following definition, as well as a helpful lemma from~\cite{HJ}.

\begin{defn}
Fix $D\ge 8$ and a plane graph $G$ with $\Delta(G)\le D$. 
A set $S\subseteq V(G)$ is \Emph{$D$-big} if $S$ is a clique in $G^2$ and $|S|\ge
D+20$.
\end{defn}

We also need the following structural lemma~\cite{HJ}.
We discuss it a bit more at the end of the paper.

\begin{lemA}[\cite{HJ}] 
Let $G$ be a simple planar graph.  There exists a vertex $v$ with $s$ neighbors
$w_1,\ldots,w_s$ such that $d(w_1)\le\cdots\le d(w_s)$ and 
one of the following holds:
\begin{enumerate}
\item $s\le 2$;
\item $s= 3$ with $d(w_1)\le 10$;
\item $s= 4$ with $d(w_1) + d(w_2)\le 15$ and $d(w_2)\le10$;
\item $s= 5$ with $d(w_1)+d(w_2)+d(w_3)\le 18$ and $d(w_3)\le 7$.
\end{enumerate}
\end{lemA}

\section{Proofs}
In the proof of our next lemma, we will often use the following easy observation.

\begin{rem}
\label{rem1}
Let $\C_1$ and $\C_2$ be cycles in a plane graph, possibly with some vertices
in common.  Fix $v_1,v_2\in V(G)$ such that $v_1$ and $v_2$ lie in distinct
components of $G-V(\C_i)$ for each $i\in\{1,2\}$.  If $\dist(v_1,v_2)\le 2$,
then $v_1$ and $v_2$ have a common neighbor in $V(\C_1)\cap V(\C_2)$.
\end{rem}
\begin{proof}
Since $v_1$ and $v_2$ lie in distinct components of $G-V(\C_1)$, 
we have $\dist(v_1,v_2)\ge 2$.
Since $\dist(v_1,v_2)\le 2$, clearly $v_1$ and $v_2$ have a common neighbor;
call it $w$.  For each $j\in\{1,2\}$, if $w\notin V(\C_j)$, then $\C_j$
separates $w$ from either $v_1$ or $v_2$; so $w$ is not a common neighbor of
$v_1$ and $v_2$, a contradiction.
\end{proof}

Most of the work proving the Big Clique Lemma goes into proving the following
more technical result.

\begin{lem}
\label{lem1}
Fix an integer $D\ge 19$, a plane graph $G$ with $\Delta(G)\le D$, and
$S\subseteq V(G)$ such that $S$ is a maximal $D$-big set.
If $d(v)+d(w)\ge D+3$ for all $v,w\in V(G)$ such that $vw\in E(G)$ and $v\notin
S$ and $d(v)\le 5$, then $G$ contains vertices $x,y,z$ such that $S=\{v\in V(G):
|N[v]\cap\{x,y,z\}|\ge 2\}$.
\end{lem}

Before proving Lemma~\ref{lem1}, we provide some intuition.  Lemma~A immediately
implies that $\chi(G^2)\le 2\Delta(G)+19$, whenever $G$ is a plane graph with
$\Delta(G)\ge 13$.  To see this, let $v$ be a $5^-$-vertex guaranteed by the
lemma, and form $G'$ from $G$ by contracting an edge $vw$, where $d(w)\le 10$ if
$d(v)\ge 3$.  By induction we color $G'^2$, and extend the coloring to $v$ using
at most $2\Delta(G)+19$ colors.  Of course, this also implies that
$\omega(G^2)\le 2\Delta(G)+19$.  So our plan is to refine this argument.  Given
any $D$-big set $S$ in $G$, we find a vertex $v$, as in Lemma~A, that lies in
$S$.  (It is our need to have $v\in S$ that leads to the hypothesis $d(v)+d(w)\ge
D+3$ for various edges $vw$.  After proving Lemma~\ref{lem1}, we show how to
remove this hypothesis; essentially, we repeatedly contract away edges that
violate it, and apply Lemma~\ref{lem1} to the resulting graph.)
Let $x$ and $y$ be two neighbors of $v$ with largest degrees.  We show
that $S\cap((N(x)\setminus N(y))\cup(N(y)\setminus N(x)))$ must have order at
most $\Delta(G)$; otherwise, not all vertices of $S$ can be pairwise adjacent in 
$G^2$, a contradiction.

\begin{proof}[Proof of Lemma~\ref{lem1}.]
Fix $D$, $G$, and $S$ satisfying the hypotheses of the lemma.
\begin{clm}
\label{clm1}
There exist $v\in S$ and $x,y\in N(v)$ such that 
$|S\setminus (N[x]\cup N[y])|\le 11$.
\end{clm}
\begin{clmproof}
Form \Emph{$G'$} from $G$ by repeating the following step as long as possible: If any
$10^-$-vertex $w$ lies on a $4^+$-face $f$, then add a chord of that face incident
to $w$ (as long as that chord is not already present, but embedded outside $f$).
Let $v$ be the $5^-$-vertex guaranteed in $G'$ by Lemma~A,\aside{$v$, $x$, $y$}
and choose $x,y\in N(v)$ to be two neighbors of $v$ with largest degrees.  In
particular, all other neighbors of $v$ are $10^-$-vertices and $\sum_{w\in
N(v)\setminus\{x,y\}}d(w)\le 18$.  If $v\notin S$, then $d_G(v)\ge 3$, since by
hypothesis $d_G(v)+d_G(w)\ge D+3$, for each edge $vw\in E(G)$.  Since
$d_G(v)\ge 3$, there exists $w\in N_G(v) \setminus \{x,y\}$.  By Lemma~A, we
have $d_G(w)\le d_{G'}(w)\le 10$.  Thus, $D+3\le d_G(v)+d_G(w)\le 10+5$, 
contradicting $D\ge 19$.  Thus, $v\in S$.

Since $v\in S$, we have $S\subseteq N_{G^2}[v]$.  By Lemma~A, we
immediately have $|S\setminus(N[x]\cup N[y])|\le 18$.  But we strengthen this
bound as follows.  For each $w\in N(v)\setminus\{x,y\}$, the graph $G'$ contains edges
from $w$ to at least two other neighbors of $v$ (the predecessor and successor of $w$ in
the cyclic order around $v$).  If not, then the edge $vw$ would lie on a
$4^+$-face $f$, and we would have added a chord of $f$ (incident to $v$ or $w$) when forming $G'$.
If $d(v)\le 4$, then $|S\setminus(N[x]\cup N[y])|\le \max\{10-2,15-2(2)\}=11$, so
we are done.  Assume instead that $d(v)=5$.  Now at least two $10^-$-neighbors
of $v$ are successive in the cyclic order, so also have a common neighbor other
than $v$; again, if not, then they would lie on a common $4^+$-face $f$, and we
would have added a chord of $f$.  Thus,  $|S\setminus(N[x]\cup N[y])|\le
18-3(2)-1=11$.
\end{clmproof}

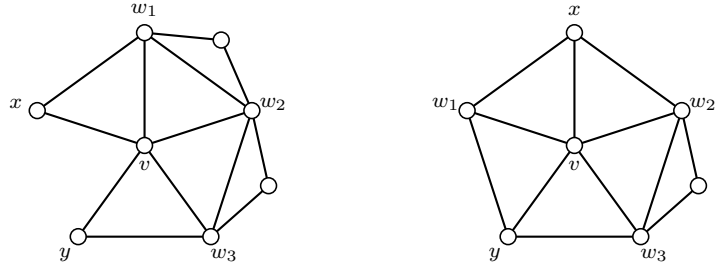
\begin{figure}[!h]
\centering
\begin{tikzpicture}[thick, scale=.6]
\tikzstyle{uStyle}=[shape = circle, minimum size = 6pt, inner sep = 0pt,
outer sep = 0pt, draw, fill=white, semithick]
\tikzstyle{sStyle}=[shape = rectangle, minimum size = 4.5pt, inner sep = 0pt,
outer sep = 0pt, draw, fill=white, semithick]
\tikzstyle{lStyle}=[shape = circle, minimum size = 4.5pt, inner sep = 0pt,
outer sep = 0pt, draw=none, fill=none]
\tikzset{every node/.style=uStyle}
\def\rad{2.5cm}
\def\off{.45cm}

\draw (0,0) node (v) {};
\draw (v) ++ (0,-\off) node[lStyle] {\footnotesize{$v$}};

\foreach \ang/\name in {18/{w_2}, 90/{w_1}, 162/x, 234/y, 306/{w_3}}
\draw (\ang:\rad) node (\name) {} (\ang:1.2*\rad) node[lStyle]
{\footnotesize{$\name$}};

\foreach \to/\from in {v/{w_1}, v/{w_2}, v/{w_3}, v/x, v/y, {w_1}/{w_2},
{w_2}/{w_3}, {w_1}/x, y/{w_3}}
\draw (\from) -- (\to);

\draw node (x_1) at (barycentric cs:{w_1}=1,{w_2}=1,v=-.6) {};
\draw node (x_2) at (barycentric cs:{w_2}=1,{w_3}=1,v=-.6) {};
\draw (w_1) -- (x_1) -- (w_2) -- (x_2) -- (w_3);

\begin{scope}[xshift=3.75in]

\draw (0,0) node (v) {};
\draw (v) ++ (0,-\off) node[lStyle] {\footnotesize{$v$}};

\foreach \ang/\name in {18/{w_2}, 90/x, 162/{w_1}, 234/y, 306/{w_3}}
\draw (\ang:\rad) node (\name) {} (\ang:1.2*\rad) node[lStyle]
{\footnotesize{$\name$}};

\foreach \to/\from in {v/{w_1}, v/{w_2}, v/{w_3}, v/x, v/y, x/{w_2},
{w_2}/{w_3}, {w_1}/x, y/{w_3}}
\draw (\from) -- (\to);

\draw node (x_2) at (barycentric cs:{w_2}=1,{w_3}=1,v=-.6) {};
\draw (w_1) -- (y) (w_2) -- (x_2) -- (w_3);
\end{scope}

\end{tikzpicture}
\caption{Two examples of $N[v]$ in Claim~\ref{clm1} and some edges that must
appear on faces incident to these vertices, due to the diagonalizing step
before we apply Lemma~A.\label{clm1-fig}}
\end{figure}

Let $X:=S\cap (N(x)\setminus N[y])$ and $Y:=S\cap (N(y)\setminus N[x])$ and
$Z:=S\cap N(x)\cap N(y)$ and $S_{xy}:=S\cap\{x,y\}$.\aside{$X, Y, Z, S_{xy}$}

\begin{clm}
\label{clm2}
$|X|\ge 7$ and $|Y|\ge 7$.
Either (i) $S\subseteq N[x]\cup N[y]$ or else
(ii) $|X|+|Y|+|S_{xy}|\ge |S|-15$.  
\end{clm}
\begin{clmproof}
By definition $S= X\udot Y\udot Z\udot S_{xy}\udot (S\setminus(N[x]\cup N[y]))$;
here $\udot$ denotes disjoint union.
By Claim~\ref{clm1}, we have
$|S\setminus (N[x]\cup N[y])|\le 11$.
Thus, $|X|=|S|-(|Y|+|Z|)-|S_{xy}|-|S\setminus(N[x]\cup N[y])|\ge
|S|-D-2-11 \ge (D+20)-(D+13) = 7$.
Similarly, $|Y|\ge 7$.
This proves the first statement.

If $|Z|\le 4$, then $|X|+|Y|+|S_{xy}| = |S|-|Z|-|S\setminus (N[x]\cup N[y])|
\ge |S| - 4 - 11$. 
Now instead assume $|Z|\ge 5$; see the left of Figure~\ref{clms2,4-fig}.  
Assume also $S\nsubseteq N[x]\cup N[y]$; otherwise, we are done.  
Denote the vertices of $Z$ by $z_1,z_2,\ldots$ in
cyclic order around $x$ and $y$.  For each $i\in[|Z|]$, let $R_i$ (short for
\Emph{region} $i$) denote the portion of the plane bounded by 
the 4-cycle $xz_iyz_{i-1}$ that contains no other vertex $z_j$; here subscripts
are modulo $|Z|$.
Since $S\nsubseteq N[x]\cup N[y]$, there exists $w\in S\setminus
(N[x]\cup N[y])$.  
By symmetry, assume $w$ lies in $R_1$.  Now 
$\dist(w,z_3)\ge 3$, by Remark~\ref{rem1};
since $w,z_3\in S$, this contradicts that $S$ \mbox{is a clique in $G^2$.}
\end{clmproof}

Denote the vertices of $X$ and $Y$, respectively, by
$x_1,\ldots,x_k$ and $y_1,\ldots,y_{\ell}$;\aside{$x_i$, $y_j$} assume these appear 
in cyclic (clockwise) order around $x$ and cyclic (counterclockwise) order around $y$.  
Further, we assume that $v$ precedes both $x_1$ and $y_1$ in these orders.
Figure~\ref{clm5-fig} shows an example.

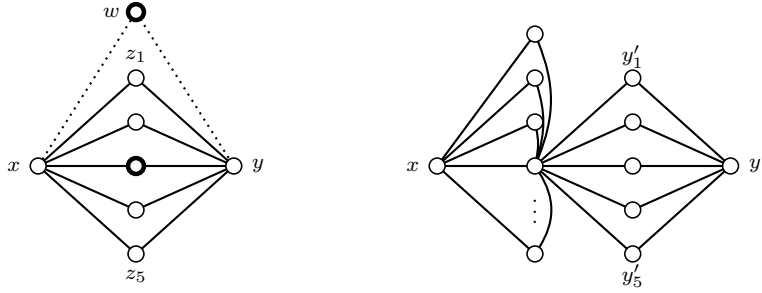
\begin{figure}[!h]
\centering
\begin{tikzpicture}[thick, scale=.65, yscale=.9]
\tikzstyle{uStyle}=[shape = circle, minimum size = 6pt, inner sep = 0pt,
outer sep = 0pt, draw, fill=white, semithick]
\tikzstyle{sStyle}=[shape = rectangle, minimum size = 4.5pt, inner sep = 0pt,
outer sep = 0pt, draw, fill=white, semithick]
\tikzstyle{lStyle}=[shape = circle, minimum size = 4.5pt, inner sep = 0pt,
outer sep = 0pt, draw=none, fill=none]
\tikzset{every node/.style=uStyle}
\def\rad{2.5cm}
\def\off{.5cm}

\foreach \i in {1,...,5}
\draw (0,3-\i) node (z\i) {};

\draw (-2,0) node (x) {} (2,0) node (y) {};

\foreach \i in {1,...,5}
\draw (x) -- (z\i) -- (y);

\draw (x) ++ (-\off,0) node[lStyle] {\footnotesize{$x$}};
\draw (y) ++ (\off,0) node[lStyle] {\footnotesize{$y$}};
\draw (z1) ++ (0,\off) node[lStyle] {\footnotesize{$z_1$}};
\draw (z5) ++ (0,-\off) node[lStyle] {\footnotesize{$z_5$}};
\draw (0,3.5) node[line width=.6mm] (w) {} ++ (-\off,0) node[lStyle] {\footnotesize{$w$}};
\draw (z3) node[line width=.6mm] {};
\draw[dotted] (x) -- (w) -- (y);

\begin{scope}[xshift=4in]

\foreach \i in {1,...,5}
\draw (0,3-\i) node (y\i) {};

\draw  (-4,0) node (x) {} -- (-2,0) node (x4) {} (2,0) node (y) {};

\foreach \i in {1,...,5}
\draw (x4) -- (y\i) -- (y);

\foreach \i in {1,2,3,6}
\draw (-2,4-\i) node (x\i) {} -- (x);

\draw (x4) edge[bend right=20] (x1);
\draw (x4) edge[bend right=15] (x2);
\draw (x4) edge[bend right=10] (x3);
\draw (-2,-.85) node[lStyle] {\footnotesize{$\vdots$}};
\draw (x4) edge[bend left] (x6);

\draw (x) ++ (-\off,0) node[lStyle] {\footnotesize{$x$}};
\draw (y) ++ (\off,0) node[lStyle] {\footnotesize{$y$}};
\draw (y1) ++ (0,\off) node[lStyle] {\footnotesize{$y'_1$}};
\draw (y5) ++ (0,-\off) node[lStyle] {\footnotesize{$y'_5$}};
\end{scope}
\end{tikzpicture}
\caption{Left: The case $S\nsubseteq N[x]\cup N[y]$ and $|Z|\ge 5$ in the proof of the 
second statement in Claim~\ref{clm2}.  Right: The proof of Claim~\ref{clm3}, where we suppose, to the contrary,
that there exists $x_i$ such that $|N(x_i)\cap Y|\ge 5$.%
\label{clms2,4-fig}}
\end{figure}

\begin{clm}
If there exists $z$ such that $X\cup Y\subseteq N[z]$, then the lemma holds.
\label{clm2.5}
\end{clm}
\begin{clmproof}
Assume such a \Emph{$z$} exists.  
By Claim~\ref{clm2}, either (i) $S\subseteq N[x]\cup N[y]$ or else (ii)
$|X|+|Y|+|S_{xy}|\ge |S|-15$.  
In the latter case, $d(z)\ge |X|+|Y|-|\{z\}|\ge |S|-15-|S_{xy}|-1\ge |S| - 18 \ge
D+2 > D$, a contradiction.  
So assume instead that $S\subseteq N[x]\cup N[y]$.
We will show that $S=\{w:|N[w]\cap\{x,y,z\}|\ge 2\}$.

We first show that $S\subseteq\{w:|N[w]\cap\{x,y,z\}|\ge 2\}$.
Since $S\subseteq N[x]\cup N[y]$, we have
$S= X\udot Y\udot Z \udot S_{xy} \subseteq (N(x)\cap N[z])\cup (N(y)\cap N[z])\cup
(N(x)\cap N(y))\cup S_{xy}\subseteq \{w:|N[w]\cap\{x,y,z\}|\ge 2\}\cup S_{xy}$.
If $x\notin N(y)\cup N(z)$, then $x\notin S$ as follows.
Since $|Y|\ge 7$ and $Y\subseteq N(y)\cap N[z]$, there exists distinct vertices
$y'_1,\ldots,y'_5\in N(y) \cap N(z)$.  Thus, there exists
$y'_j\in Y$ such
that $\dist(x,y_j)\ge 3$.  The proof of this fact mirrors the final paragraph
proving Claim~\ref{clm2}.  By definition, $y'_j\in Y\subseteq S$. 
So $\dist(x,y'_j)\ge 3$ implies that $x\notin S$.  Similarly, if $y\notin
N(x)\cup N(z)$, then $y\notin S$.  Hence, as desired,
$S\subseteq\{w:|N[w]\cap\{x,y,z\}|\ge 2\}$.

We now show that $\{w:|N[w]\cap\{x,y,z\}|\ge 2\}\subseteq S$.
As proved above, if $w\in S$, then $|N[w]\cap \{x,y,z\}|\ge 2$.
So fix arbitrary vertices $w_1,w_2$ such that $|N[w_i]\cap\{x,y,z\}|\ge 2$ for
each $i\in\{1,2\}$.  By Pigeonhole, $N[w_1]\cap N[w_2]\ne\emptyset$.  So
$w_1w_2\in E(G^2)$.  Since $S$ is a maximal clique in $G^2$, we have
$S=\{w:|N[w]\cap\{x,y,z\}|\ge 2\}$.
\end{clmproof}

\begin{clm}
\label{clm3}
For all $x_i\in X$ and $y_j\in Y$ we have $|N(x_i)\cap Y|\le 4$ and $|N(y_j)\cap X|\le 4$.
\end{clm}
\begin{clmproof}
Suppose, to the contrary, that there exists $x_i\in X$ with $|N(x_i)\cap Y|\ge 5$.
Fix $y_1',\ldots,y_5'\in N(x_i)\cap Y$.  By symmetry, we assume that
$y_1',\ldots,y_5'$ are in order cyclically around $y$; see the right of
Figure~\ref{clms2,4-fig}.
Now $x_iy_2'yy_4'$ separates $y_3'$ from $y_1'$, $y_5'$, and all vertices in $X\setminus \{x_i\}$.
Similarly, $x_iy_1'yy_5'$ separates $y_2'$, $y_3'$, and $y_4'$ from all vertices in $X\setminus \{x_i\}$.
Thus, $X\subseteq N[x_i]$ by Remark~\ref{rem1}.
If $i\ge 4$, then (a) $x_{i-1}$ is separated from $Y$ by both $xx_ix_{i-2}$ and
$xx_ix_{i-3}$.  And if $i\le 3$, then (b) $x_{i+1}$ is separated from $Y$ by
both $xx_ix_{i+2}$ and $xx_ix_{i+3}$ (since $|X|\ge 7$, by Claim~\ref{clm2}).
Now $Y\subseteq N[x_i]$ by Remark~\ref{rem1}.
Hence, $X\cup Y\subseteq N[x_i]$.  So, we are done by Claim~\ref{clm2.5}.
\end{clmproof}

A \Emph{chord} in $X$ is an edge $x_ix_j$ with $|i-j|>1$. 
Similarly, a \emph{chord} in $Y$ is an edge $y_iy_j$ with $|i-j|>1$. 

\begin{clm}
\label{clm4}
$G$ does not contain both a chord $x_ix_j$ in $X$ and a chord $y_{i'}y_{j'}$ in $Y$.
Thus, by symmetry, we may assume that $G$ does not contain a chord in $Y$.
\end{clm}
\begin{clmproof}
Suppose the contrary; by symmetry, assume $i<j$ and $i'<j'$.
Now $xx_ix_j$ and $yy_{i'}y_{j'}$ each separate $x_{i+1}$
from $y_{i'+1}$. So Remark~\ref{rem1} implies $\dist(x_{i+1},y_{i'+1})\ge 3$,
contradicting that $x_{i+1},y_{i'+1}\in S$.
\end{clmproof}

\begin{clm}
\label{clm5}
$G$ does not contain an edge $x_iy_j$ with $2\le i\le k-1$ and $3\le j\le \ell-2$.
\end{clm}
\begin{clmproof}
Suppose, to the contrary, that such an edge exists; see Figure~\ref{clm5-fig}.
Now $x_1$ is separated from $y_{\ell}$ by the cycle $xvyy_jx_i$.
Since $\dist(x_1,y_{\ell})\le 2$, some vertex of this cycle must be a common
neighbor $z$ of $x_1$ and $y_{\ell}$.  Clearly, $z\notin\{x,y\}$.
Recall that $G$ does not contain a chord of $Y$, by Claim~\ref{clm4}.
Thus, $z\ne y_j$.
So $z\in\{v,x_i\}$.
We can repeat the argument above for each $x_h$ with $1\le h\le i-1$.
We can also repeat it for each $x_h$ with $i+1\le h\le k$, using vertex $y_1$ in
place of $y_{\ell}$.  We can also repeat the argument for $x_1$ and $y_h$
whenever, $j+2\le h\le \ell$; and we can repeat it for $x_k$ and $y_h$
whenever $1\le h\le j-2$.  Thus, $X\cup(Y\setminus\{y_{j-1},y_{j+1}\})\subseteq
N[x_i]\cup N(v)$.
Below we will strengthen this to show that $X\cup Y\subseteq N(v)$, and finish
by Claim~\ref{clm2.5}.

Claim~\ref{clm2} implies that either (i) $S\subseteq N[x]\cup N[y]$ or else (ii)
$|X|+|Y|\ge |S|-|S_{xy}|-15\ge (D+20)-17 = D+3\ge 21$, since $D\ge 18$.
Assume (i) holds; so $S=X\udot Y\udot Z\udot S_{xy}$.  
Now $|X|+|Z|\le D$ and $|Y|+|Z|\le D$ and $|S_{xy}|\le 2$.
We sum these three inequalities, and we subtract the sum from
$2(|X|+|Y|+|Z|+|S_{xy}|)=2|S|$; this gives
$|X|+|Y|+|S_{xy}|\ge 2|S|-2D-2\ge 2(D+20)-2D-2\ge 38$.
Thus $|X|+|Y|\ge 38-|S_{xy}|\ge 36$. 
So in both case (i) and case (ii), we have $|X|+|Y|\ge \min\{36,21\}$.
Thus, either $|X|\ge 10$ or $|Y|\ge 12$.

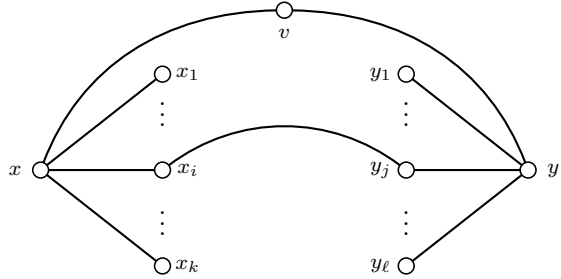
\begin{figure}[!h]
\centering
\begin{tikzpicture}[thick, scale=.675, yscale=.9]
\tikzstyle{uStyle}=[shape = circle, minimum size = 6pt, inner sep = 0pt,
outer sep = 0pt, draw, fill=white, semithick]
\tikzstyle{sStyle}=[shape = rectangle, minimum size = 4.5pt, inner sep = 0pt,
outer sep = 0pt, draw, fill=white, semithick]
\tikzstyle{lStyle}=[shape = circle, minimum size = 4.5pt, inner sep = 0pt,
outer sep = 0pt, draw=none, fill=none]
\tikzset{every node/.style=uStyle}
\def\rad{2.5cm}
\def\xscale{2.4}
\def\yscale{0.7}
\def\off{5mm}

\foreach \name/\x/\y/\flip in {{x_1}/-1/3/1, {x_i}/-1/0/1, {x_k}/-1/-3/1,
{y_1}/1/3/-1, {y_j}/1/0/-1,
{y_l}/1/-3/-1, y/2/0/1, x/-2/0/-1, v/0/5/1}
\draw (\x*\xscale,\y*\yscale) node (\name) {} ++ (\flip*\off,0) node[lStyle] {\footnotesize{$\name$}};

\draw (v) ++ (\off,0) node[draw=none, fill=white] {};
\draw (y_l) ++ (-\off,0) node[draw=none, fill=white, minimum size=9 pt] {};
\draw (y_l) ++ (-\off,0) node[lStyle] {\footnotesize{$y_{\ell}$}};
\draw (v) ++ (0,-\off) node[lStyle] {\footnotesize{$v$}};

\foreach \a/\b in {1/1, 1/-1, -1/1, -1/-1}
\draw (\a*\xscale,\b*1.7*\yscale+.2) node[lStyle] {{$\vdots$}};

\draw (x_1) -- (x) -- (x_i) (x) -- (x_k) (y_1) -- (y) -- (y_j) (y) -- (y_l);
\draw (x) edge[bend left=35] (v);
\draw (y) edge[bend right=35] (v);
\draw (x_i) edge[bend left=40] (y_j);

\end{tikzpicture}
\caption{An edge $x_iy_j$ with $2\le i\le k-1$ and $3\le j\le
\ell-2$, as in the proof of Claim~\ref{clm5}.\label{clm5-fig}}
\end{figure}

\textbf{Case 1: $|X|\ge 10$.}  
Suppose that $x_hx_i\in E(G)$ for at least two values of $h$ with $1\le h\le
i-2$.  Now the cycle $xx_hx_i$ separates $x_{i-1}$ from $Y$ for each $h$.
Remark~\ref{rem1} gives
$x_iy_h\in E(G)$ for all $y_h\in Y$; that is, $Y\subseteq N(x_i)$. 
But this contradicts Claim~\ref{clm3}, since $|Y|\ge 7$ by Claim~\ref{clm2}.
Similarly, we get a contradiction if $x_hx_i\in E(G)$ for at least two values of
$h$ with $i+2\le h\le k$. 

Recall from above that $X\subseteq N[x_i]\cup N(v)$.
If $i\ge 6$, then there exist $h', h'', h'''$ such that $1\le
h'<h''<h'''\le 4$ and $x_{h'},x_{h''},x_{h'''}\in N(v)$.  So $x_{h'}$ is
separated from $Y$ by cycles $xvx_{h''}$ and
$xvx_{h'''}$.  Now Remark~\ref{rem1} implies that $Y\subseteq N(v)$.  If instead $i\le 5$, then
there exist $h', h'', h'''$ such that $k\ge h'>h''>h'''\ge 7$ and
$x_{h'},x_{h''},x_{h'''}\in N(v)$.
Again $x_{h'}$ is separated from $Y$ by cycles $xvx_{h''}$ and $xvx_{h'''}$.
So again Remark~\ref{rem1} implies that $Y\subseteq N(v)$.  Since $Y\subseteq N(v)$, either (a)
$y_1$ is separated from $X$ by cycles $yvy_2$ and $yvy_3$ (if $j\ge 4$) or else
(b) $y_6$ is separated from $X$ by cycles $yvy_5$ and $yvy_4$ (recall from
Claim~\ref{clm2} that $|Y|\ge 7$).  In each case, $X\subseteq N(v)$.  
Thus, $X\cup Y\subseteq N(v)$, so we are done by Claim~\ref{clm2.5}.

\textbf{Case 2: $|Y|\ge 12$.}  
Suppose that $y_hx_i\in E(G)$ for at least two values of $h$ with $1\le h\le
j-3$; say $h'$ and $h''$.  Now the cycle $yy_hx_iy_j$ separates $y_{j-2}$ from
$X$ for each $h\in\{h',h''\}$.  By Remark~\ref{rem1}, each vertex of $X$ has a
common neighbor $z$ with $y_{j-2}$ among $\{y,x_i,y_j\}$.  Clearly, $z\ne y$,
and also $z\ne y_j$, since $G$ has no chords in $Y$.  Thus, $z=x_i$.  That is,
$X\subseteq N[x_i]$.  If $i\ge 4$, then
$x_3$ is separated from $Y$ by cycles $xx_ix_1$ and $xx_ix_2$.  This implies
that all of $Y$ is adjacent to $x_i$.  So $X\cup Y\subseteq N[x_i]$, and we are
done by Claim~\ref{clm2.5}.
Thus, $y_hx_i\in E(G)$ for at most one
value of $h$ with $1\le h\le j-3$.  A nearly identical argument gives that
$y_hx_i\in E(G)$ for at most one value of $h$ with $j+3\le h\le \ell$.  

Recall that $Y\setminus\{y_{j-1},y_{j+1}\}\subseteq N(x_i)\cup N(v)$.
If $j\ge
7$, then there exist $h',h'',h'''$ with $1\le h'<h''<h'''\le 4$ such
that $y_{h'},y_{h''},y_{h'''}\in N(v)$.  Thus,
$y_{h'}$ is separated from $X$ by cycles $yy_hv$ for each
$h\in\{h'',h'''\}$.  As above, we conclude that $X\subseteq N(v)$.  But now, since
$|X|\ge 7$, either $x_1$ is separated from $Y$ by cycles $xx_2v$ and $xx_3v$ or
else $x_6$ is separted from $Y$ by cycles $xx_5v$ and $xx_4v$.  In each case we
conclude that $Y\subseteq N(v)$.  Thus $X\cup Y\subseteq N(v)$, so we are done
by Claim~\ref{clm2.5}.
So assume instead that $j\le 6$.  Now the argument is essentially the same,
except that we pick $h',h'',h'''$ with $\ell\ge h'>h''>h'''\ge 9$
and $y_{h'},y_{h''},y_{h'''}\in N(v)$.
So $y_{h'}$ is separated from $X$ by cycles $yy_hv$ for each
$h\in\{h'',h'''\}$.  As above, this gives $X\subseteq N(v)$ and,
eventually, $Y\subseteq N(v)$.
Thus, $X\cup Y\subseteq N(v)$, so we are done by Claim~\ref{clm2.5}.
\end{clmproof}

\begin{clm}
There exists $z\in (N(x_4)\cap N(y_4))\setminus(X\cup Y)$.
\end{clm}
\begin{clmproof}
Recall that $|X|\ge 7$ and $|Y|\ge 7$, by Claim~\ref{clm2}.
By Claim~\ref{clm5}, we have $x_iy_j\notin E(G)$ for all $i\in \{2,\ldots,
k-1\}$ and $j\in\{3,4,5\}$.  By Claim~\ref{clm2}, we have $y_4y_j\notin E(G)$ for all $j\in
[\ell]\setminus\{3,4,5\}$.  Thus, $x_4$ and $y_4$ have a common neighbor $z$
and either (a) $z\in\{x_1,x_k\}$ or (b) $z\notin X\cup Y$ (recall that
$v\notin X\cup Y$).
Suppose that $x_1\in N(x_4)\cap N(y_4)$.  Now $xx_1x_4$ separates $x_2$ and $x_3$ 
from $y_1$.  Further, $yy_4x_1xv$ separates $x_2$ and $x_3$ from $y_1$.
Thus, $x_2,x_3\in N(x_1)$ by Remark~\ref{rem1}.  But now $xx_1x_3$ and $xx_1x_4$ both separate $x_2$
from $Y$.  Thus, $Y\subseteq N(x_1)$ by Remark~\ref{rem1}.  As a result, $y_2$ is
separated from $X\setminus \{x_1\}$ by both $yy_3x_1y_1$ and $yy_4x_1xv$.
So $X\subseteq N(x_1)$, by Remark~\ref{rem1}.  Thus, $X\cup Y\subseteq N[x_1]$
and we are done by Claim~\ref{clm2.5}.
If $x_k\in N(x_4)\cap N(y_4)$, then
the proof is nearly identical, by symmetry. 
Hence, (a) above cannot hold.  
\end{clmproof}

Suppose $|N(v)\cap X|\ge 3$.
Now $x_1$ is separated from $Y$ by $xvx_h$ for at least two values of $h$ with
$x_h\in X\setminus \{x_1\}$.  So $Y\subseteq N(v)$.
But now we can repeat this argument to show that $vyy_h$ separates $y_1$ from
$X$ for at least two values of $h$ with $y_h\in Y\setminus\{y_1\}$.  So
$X\subseteq N(v)$.  Thus, $X\cup Y\subseteq N(v)$, and we are done by
Claim~\ref{clm2.5}.
As a result, $|N(v)\cap X|\le 2$.  
Each vertex $x_i$ with $1\le i\le 3$ is separated from $y_6$ by the cycle
$xx_4zy_4yv$.  (If $v=z$, then $yy_4z$ separates $X$ either from
$\{y_1,y_2,y_3\}$ or from $\{y_5,y_6,y_7\}$; by symmetry, assume the former.
Now $y_2y_4\notin E(G)$, so $y_2v\in E(G)$.  Thus, $yy_2v$ and $yy_4v$ both
separate $y_1$ from $X$; hence $X\subseteq N[v]$, contradicting $|N(v)\cap X|\le
2$, which we proved above.)

For each $h\in\{2,3\}$, vertex $x_h$ is separated from $y_5$ by $vxx_4zy_4y$; so
$x_h$ and $y_5$ must have a common neighbor $w$ on this cycle.  Clearly,
$w\notin\{x,y\}$.  And Claim~\ref{clm5} implies $w\notin\{x_4,y_4\}$.  Thus,
$w\in \{v,z\}$.  A symmetric argument applies to $y_3$ and $x_h$ for each
$h\in\{5,6\}$.  Hence, $\{x_2,\ldots,x_6\}\subseteq N(v)\cup N(z)$.
Since $x_2\in N(v)\cup N(z)$, we know $x_1x_4\notin E(G)$.  So $x_1$ and $y_6$
have a common neighbor in $\{v,z\}$.  That is, $x_1\in N(v)\cup N(z)$.  By
symmetry, also $x_7\in N(v)\cup N(z)$.
Thus, $X\subseteq N(v)\cup N(z)$.  Since $|X\cap N(v)|\le 2$, we get that
$|X\cap N(z)|\ge 5$.  

Let $x'_1,\ldots,x'_5$ denote 5 vertices in $X\cap N(z)$, in cyclic order
around both $x$ and $z$.  But now the cycles $x'_1zx'_5x$ and $x'_2zx'_4x$ both
separate $Y$ from $x'_3$.  So $Y\subseteq N(z)$.  By swapping the roles of $X$
and $Y$, we can use the same argument to show that $X\subseteq N(z)$.  But now
$X\cup Y\subseteq N(z)$, so we are done by Claim~\ref{clm2.5}.
\end{proof}

Now we can prove the Big Clique Lemma.  For easy reference, we restate it here.

\begin{clique-lem}
Fix an integer $D\ge 19$ and a plane graph $G$ with $\Delta(G)\le D$. If
$S\subseteq V(G)$ such that $S$ is a maximal clique in $G^2$ and $|S|\ge
D+20$,
then $G$ contains vertices $x,y,z$ such that $S=\{v\in V(G):
|N[v]\cap\{x,y,z\}|\ge 2\}$.
\end{clique-lem}

\begin{figure}[!t]
\centering
\begin{tikzpicture}[thick, scale=.75, yscale=.985]
\tikzstyle{uStyle}=[shape = circle, minimum size = 6pt, inner sep = 0pt,
outer sep = 0pt, draw, fill=white, semithick]
\tikzstyle{sStyle}=[shape = rectangle, minimum size = 4.5pt, inner sep = 0pt,
outer sep = 0pt, draw, fill=white, semithick]
\tikzstyle{lStyle}=[shape = circle, minimum size = 4.5pt, inner sep = 0pt,
outer sep = 0pt, draw=none, fill=none]
\tikzset{every node/.style=uStyle}
\def\rad{2.5cm}
\def\off{.6cm}

\foreach \ang/\name/\newname in {90/z/x, 210/y/z, 330/x/y}
\draw (\ang:\rad) node[uStyle] (\name) {} ++ (\ang:\off) node[lStyle]
(\newname1) {\footnotesize{$\newname$}};
\draw (z1) ++ (-.8*\off,0) node[lStyle] {\footnotesize{$w=~~$}};

\foreach \ang/\name in {150/X, 270/Z, 30/Y}
{
\draw[rotate=\ang] (0:.75*\rad) ellipse (19pt and 10pt) {}; 
\draw (\ang:1.2*\rad) node[lStyle] {}; 
}

\foreach \name/\ang in {y/153, y/267, x/273, x/27, z/33, z/147}
{
\draw (\name) -- (\ang:\rad) (\name) -- (\ang:.49*\rad);
}

\draw[white, line width=.4mm] (y) -- (267:\rad);
\draw (y) -- (261:.85*\rad);
\draw (y) node {};
\draw (y) --++ (280:1.2cm) node (v) {} -- (268:1.01*\rad);
\draw (v) -- (261:.87*\rad) (v) ++ (-.9*\off,0) node[lStyle] {\footnotesize{$v$}};
\draw (v) ++ (358:1.95cm) node (u) {} ++ (0,-.9*\off) node[lStyle]
{\footnotesize{$u$}};
\draw[rotate=270] (0:.75*\rad) ellipse (19pt and 10pt) {}; 
\draw (u) node[fill=white] {};

\draw (0,0) node[lStyle] {\footnotesize{$G$}};

\begin{scope}[xshift=-3.75in]

\foreach \ang/\name/\newname in {90/z/x, 210/y/z, 330/x/y}
\draw (\ang:\rad) node[uStyle] (\name) {} ++ (\ang:\off) node[lStyle]
(\newname1) {\footnotesize{$\newname$}};
\draw (z1) ++ (-\off,.8mm) node[lStyle] {\footnotesize{$w'=~~$}};

\foreach \ang/\name in {150/X, 270/Z, 30/Y}
{
\draw[rotate=\ang] (0:.75*\rad) ellipse (19pt and 10pt) {}; 
\draw (\ang:1.2*\rad) node[lStyle] {}; 
}

\draw (0,0) ++ (270:2.5cm) node (u) {} ++ (0,-.9*\off) node[lStyle]
{\footnotesize{$u$}};

\foreach \name/\ang in {y/153, y/267, x/273, x/27, z/33, z/147}
{
\draw (\name) -- (\ang:\rad) (\name) -- (\ang:.49*\rad);
}

\draw (0,0) node[lStyle] {\footnotesize{$G'$}};

\end{scope}
\end{tikzpicture}
\caption{Left: A clique $S'$ in $G^2$ such that $S'=\{w: |N[w]\cap \{x,y,z\}|\ge
2\}$. Right: The preimage of $S'$ in $G$.\label{final-fig}}
\end{figure}
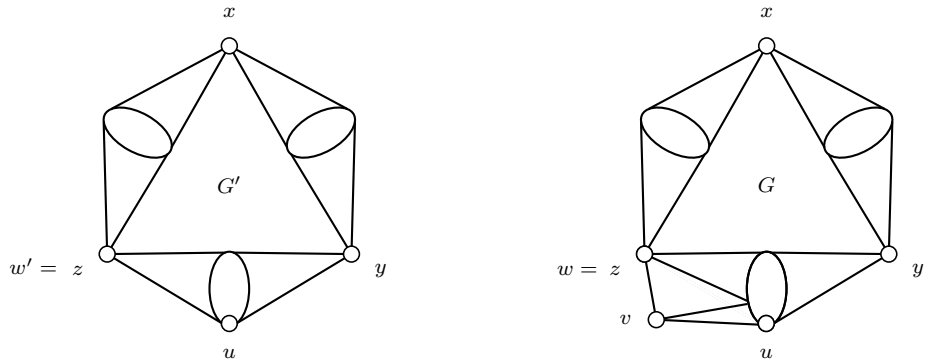

\begin{proof}
Starting from $G$, we will repeatedly contract edges until we reach a graph $G'$
such that $d_{G'}(v)+d_{G'}(w)\ge D+3$ for all $v,w\in V(G')$ such that $vw\in
E(G')$ and $v\notin S$ and $d_{G'}(v)\le 5$.  By Lemma~\ref{lem1}, the
conclusion of the present lemma holds for such a graph $G'$.  After this, we
will show that if $S$ satisfies the
desired description in $G'$, then it also does so in $G$.  Our proof is by
induction on the number $t$ of edges contracted to form $G'$ from $G$.
The base case, $t=0$, follows from Lemma~\ref{lem1}.  So assume $t\ge 1$.

By the induction hypothesis, the lemma holds for all graphs $G$ which reach a
graph satisfying the hypothesis of Lemma~\ref{lem1} after at most $t-1$ edge
contractions.  So we only need to consider the case $t=1$.  That is, $G'$ is
formed from $G$ by a single edge contraction.

Suppose there exist $v,w\in V(G)$ with $vw\in E(G)$ and $v\notin S$ and
$d(v)\le 5$ and $d(v)+d(w)\le D+2$.  Form $G'$ from $G$ by contracting edge
$vw$ (taking the underlying simple graph if this yields parallel edges);
call the resulting new vertex $w'$, and if $w\in S$ (in $G$), then say also that
$w'\in S$ in $G'$.  Since $d_G(v)+d_G(w)\le D+2$, we have $d_{G'}(w')\le D+2-2=D$.
Thus, $\Delta(G')\le \max\{\Delta(G),d_{G'}(w')\}\le D$.  And all
$u_1,u_2\in V(G)\setminus\{v\}$ satisfy $\dist_{G'}(u_1,u_2)\le \dist_G(u_1,u_2)$.
Thus, all $u_1,u_2\in S$, satisfy $\dist_{G'}(u_1,u_2)\le \dist_G(u_1,u_2)\le 2$.
So $S$ is a clique in $(G')^2$.  Let $S'$ be a maximal clique in $(G')^2$
containing $S$.

By the induction hypothesis, there exist $x,y,z\in V(G')$ such that 
$S'=\{u\in V(G'):|N_{G'}[u]\cap\{x,y,z\}|\ge 2\}$.  If $w'\notin\{x,y,z\}$, then
$N_G[u]\cap\{x,y,z\} = N_{G'}[u]\cap\{x,y,z\}$ for all $u\in V(G)\setminus\{v,w\}$.
And $N_G(w)\cap\{x,y,z\}\subseteq N_{G'}(w')\cap\{x,y,z\}$.
Thus, $S \subseteq \{u\in V(G):|N_G[u]\cap\{x,y,z\}|\ge 2\}$.  And, by the
maximality of our choice of $S$, this set containment holds with equality.
So we are done.

So assume instead that $w'\in\{x,y,z\}$.  
By symmetry, say $w'=z$; see Figure~\ref{final-fig}.
Suppose $u\in S'$ and $|N_G[u]\cap\{x,y,z\}|\le 1$.  We must show that $u\notin S$.
Since $u\in S'$, we must have $u\in N_G(v)\setminus\{w\}$ and $u\notin N_G(w)$.
By symmetry, we assume $u\notin N_G(x)\cup N_G(z)$.
Note that $|N_G(x)\cap N_G(z)\cap S|\ge |S|-|N_G(y)\cap
S|-|\{x,y,z\}|-(d_G(v)-1) \ge D+20 - D - 3 - 4 =13$.
Let $w_1,\ldots,w_5$ denote 5 vertices in $N_G(x)\cap N_G(w)\cap S$, in cyclic
order around $w$ (and around $x$), such that $v$ lies in the ``exterior'' of the
cycle $ww_1xw_5$, but no other $w_i$ lies in this exterior.  Now cycles 
$ww_1xw_5$ and $ww_2xw_4$ both separate $u$ from $w_3$.  Thus,
$\dist_G(u,w_3)\ge 3$.  Since $S$ is a clique in $G^2$, this implies that
$w\notin S$, as desired.
\end{proof}

The proof of Lemma~A, from~\cite{HJ}, is nearly 9 pages long.  (In fact, the
version proved there gives extra detail about the degrees of the vertices in
$N(v)\setminus\{x,y\}$.)  However, that lemma is a refinement of various
earlier results such as~\cite[Lemma~2.1]{vdHM03} and~\cite[Theorem~2]{BKPY}, either
of which would work in place of Lemma~A, still yielding a version of the Big
Clique Lemma, but with slightly worse constants.  A particularly easy
structural lemma that could serve in this role is~\cite[Lemma~3.5]{guide},
which says that every planar graph has a $5^-$-vertex with at most two
$12^+$-neighbors.  Its proof is under half a page.

\section*{Acknowledgments}
Thanks to Havet, van den Heuvel, McDiarmid, and Reed~\cite{HvdHMR} for
suggesting the general idea of this proof.  Thanks also to Hell and
Seyffarth~\cite{HellS}, who computed the maximum order of diameter 2 planar graphs
with a given maximum degree.  Ideas from their paper were crucial for our proof.
Thanks to Reem Mahmoud for reading an earlier version of this paper.  Her
feedback helped improve our presentation.

\scriptsize{
\bibliographystyle{habbrv}
{\bibliography{GraphColoring}}
}

\end{document}